\newif\ifdraft
\draftfalse

\ifdraft
\documentclass{amsart}
 \usepackage{a4wide}
\else
\documentclass[12pt]{amsart}
 \usepackage{a4wide}
 \fi
\usepackage{mathrsfs}
\usepackage{amssymb}                    
\usepackage{eepic}

 \date{2000-07-27}

\newcommand{\I}{I}
\newcommand {\x}{{\tt x}}
\newcommand {\y}{{\tt y}}

\newcommand{\CI}{\C_\I}
\newcommand{\C}{\mathscr{C}}
\newcommand{\D}{\mathscr{D}}
\newcommand{\F}{\mathscr{F}}

\newcommand{\on}{\restriction}
\newcommand{\ran}{{\rm ran}}
\renewcommand{\O}{{\mathscr O}}
\renewcommand{\P}{{\mathscr P}}

\makeatletter
\newcommand{\KNUTHcases}[1]
{\left \{\,\vcenter {\normalbaselines \m@th
\ialign {$##\hfil $&\quad ##\hfil \crcr #1\crcr }}\right .
}  
\makeatother

\def \itm#1 {\item[(#1)]}

\newcommand{\oo}[1]{\O^{(#1)}}
\newcommand{\oc}[1]{\C^{(#1)}}

\newcommand{\oon}{\O^{\langle 1 \rangle}}
\newcommand{\nd}{\mathord{\raise1pt\hbox{$\nabla$}\!\!\Delta}}
\renewcommand{\H}{\U^{(2)}}
\newcommand{\U}{{\mathscr U}}
\renewcommand{\P}{{\mathscr P}}
\newcommand{\Q}{{\mathscr Q}}

\newcommand{\Pol}{{\rm Pol}}

\swapnumbers
\theoremstyle{remark}

\def\xx#1 {\newtheorem{#1}[thm]{#1}}
\xx Lemma
\xx {Main Lemma}
\xx Theorem
\xx Corollary
\xx Idea
\xx {First Step}
\xx Consequence
\xx Example
\xx Definition
\xx Fact
\xx {Definition and Fact}
\xx{Open Questions}
\xx Setup
\xx Notation
\xx {More Notation}
\xx Remark 
\xx Acknowledgment
\xx Conclusion

\title{Clones on regular cardinals}

\thanks{This paper is available 
from {\tt arXiv.org} and 
 from the authors' homepages}


\author{Martin Goldstern}

\address{Algebra, TU Wien
\\
Wiedner Hauptstra\3e 8-10/118.2
\\
A-1040 Wien}
\email{Martin.Goldstern@tuwien.ac.at}
\urladdr{http://info.tuwien.ac.at/goldstern/}

\thanks{The first author is grateful to  the Hebrew University
of Jerusalem for the hospitality during his visit,
and  to the Austrian Science foundation for supporting the joint research
under FWF grant P13325-MAT}

\author{Saharon Shelah}
\thanks{The second author is supported by the
   Israel Science Foundation
   founded by the Israel Academy of Sciences and Humanities}

\address{Mathematics\\ Hebrew University of Jerusalem\\ 91904
Jerusalem\\ Israel} 
\email{shelah@math.huji.ac.il}
\urladdr{http://math.rutgers.edu/\char`\~shelah}

\begin{document}

\begin{abstract}
We investigate the structure of the lattice of clones on an infinite
set~$X$.  We first observe that ultrafilters naturally induce clones;
this yields a simple proof of Rosenberg's theorem: 
 there are $2^{2^{\lambda}}$
many maximal (= ``precomplete'') clones on a set of size~$\lambda$. The
clones we construct do not 
contain all unary functions. 

We then investigate clones that do contain all unary functions.  
Using a strong negative partition theorem  from pcf theory 
we show that for many cardinals $ \lambda $
  (in particular, for all successors of regulars) 
 there are  $2^{2^\lambda }$ many such clones on a set of 
 size~$\lambda $. 

Finally, we show that on a weakly compact cardinal there are exactly 2
maximal  clones  which contain all unary functions. 
\end{abstract}

\maketitle
\section{Introduction}
\label{section.intro}

\begin{Definition}  
Let $X$ be a nonempty set.   The {\em full clone} on~$X$, called $\O$
or $\O(X)$ is the set of all finitary 
functions on~$X$:  $\O = \bigcup_{n=1}^ \infty \oo n$, 
where $\oo n$ is the set of
all functions from $X^n$ into~$X$. 
\\
A {\em clone} (on $X$)
 is a set $\C \subseteq \O $
  which contains all
projections and is closed under composition. 
That is, 
\begin{enumerate}
\item For all $1\le k\le n$, the function $\pi^n_k:X^n\to X$, 
$\pi^n_k(x_1,\ldots, x_n) = x_k$, is in~$\C$.  
\item  whenever $f_1, \ldots,
f_k\in \C\cap \oo n$, $g\in \C\cap \oo k$, then the function 
$$(x_1, \ldots, x_n) \mapsto g(f_1(x_1, \ldots, x_n), \ldots, 
f_k(x_1, \ldots, x_n)) $$
(which we sometimes call $g(f_1,\ldots, f_k)$) 
is also in~$\C$. 
\end{enumerate}

Alternatively, $\C$ is a clone if $\C$ is the set of term functions of
some universal  algebra over~$X$. 
\end{Definition}

The set of clones over $X$ forms a complete algebraic lattice with
largest element~$\O$.  The coatoms of this lattice are called
``precomplete clones'' or ``maximal clones''.

Many results for clones on finite sets, and in particular a
classification of all precomplete clones on finite sets can be found in 
\cite{Sz86}.

 Rosenberg proved in  \cite{Ro76}
 that if $X$ is an infinite set of cardinality
 $\lambda $ then there are $2^{2^\lambda }$ many precomplete clones
 on~$X$. In section \ref{section.easy} we will give a short new  proof of
 this  theorem, using ultrafilters. 

Let~$\oon$, the ``full unary clone'',  
 be the clone generated by~$\oo 1$, i.e., 
the set of functions which depend only on one argument: 
$$ 
\oon := 
\{ f \circ \pi^n_k :   f\in \oo 1 , 1 \le k \le n \}
$$

The clones that we construct in section \ref{section.easy}, as well as the
clones in the  family constructed by Rosenberg, all have the property
that they induce a maximal proper submonoid of the  monoid $ \oo 1$ 
of all
 unary
functions.  This raises the following question:  What is the structure
of those clones that contain the full monoid of all unary functions,
i.e., the interval~$[\oon , \O]$?   In particular, what can we say
about the precomplete elements in this interval? 

If $X$ is a {\bf finite} set with $k$ elements, then it is known that 
 this interval is actually a finite chain (with $k+1$ many elements).
 In particular, there is a unique precomplete clone above the full unary
 clone, namely,  the set of all functions which are either essentially
 unary or not onto. 

We now turn to {\bf infinite} sets.  Again we will be mainly
interested in the maximal or ``precomplete'' clones above $\oon$. 
Since $\O$ is finitely generated over~$\O_1$, it is clear that the
interval  $[\oon, \O]$ is dually atomic, that is, every $\C\in  
[\oon, \O)$ is contained in some precomplete $\C' 
 \in[\oon, \O)$.  
(See fact \ref{fact.zorn}.)

For the case of countable~$X$, Gavrilov  proved in \cite{Ga65}
that there are exactly
2 precomplete clones in this interval, and 
Davies and Rosenberg (see \cite{DR85})
 gave an  explicit example
of one precomplete clone in this interval for every infinite~$X$.

It turns out that (for any infinite set $X$ of regular cardinality),
the clones on $X$ above $\oo 1$
 can be naturally divided into 2 classes, depending
on whether the binary functions of the clone are all ``almost unary''
or if there is a ``heavily binary'' function among them (see
definitions \ref{def.almost} and \ref{def.heavily}).

     In section \ref{section.nobinary} we show 
that among the clones whose binary part is almost unary, there is a
unique precomplete clone (namely, the clone from \cite{DR85}).  

Finally, we discuss the case which was hitherto unknown, and which
turns out to be the most interesting from the set theoretical point of
view: clones with heavily binary functions.   The structure of the set
of these clones depends on partition properties of the cardinality of
the underlying set: 
\begin{enumerate}
\item If the cardinality of the underlying set is a weakly compact
cardinal (or ${\aleph_0}$),
 then there is a unique precomplete clone in $[\O^{\langle
1\rangle},\O]$ which is heavily binary (so altogether there are exactly two
precomplete clones above $\oo 1$)\\
 This result, which generalizes Gavrilov's theorem for ${\aleph_0} $, 
 is proved in section \ref{section.wc}.
\item If the cardinality $\lambda $ of the underlying set satisfies a
certain negative partition property $Pr(\lambda  )$ (see also \ref{def.pr})
 --- in particular,
we know $Pr(\kappa ^+)$ for all regular $\kappa $), then there are
$2^{2^\lambda  }$ many precomplete clones above $\oo 1 $ which are heavily
binary. 
\\
This result is proved in section \ref{section.many}.
\end{enumerate}
In an appendix we briefly discuss partition relations and 
the combinatorial principle
$Pr(\lambda)$.

All sections of the paper can be read independently, but they all rely
on notation, facts and concepts  established in this
introduction.

We plan to investigate clones on singular cardinals in a separate paper.

\begin{Notation} \label{notation}
We fix an infinite set~$X$.  For $n\in\{1,2,\ldots\}$ we write $\oo n$
for the set of all functions from $X^n$ to~$X$, $\O =
\bigcup_{n=1}^\infty \oo n$. 

  For any set of functions $\F \subseteq \O$
 we let $cl(\F)$ be the
smallest clone containing $\F$ as well as all unary functions. 

We will write $\lambda  = |X| $ for the cardinality of~$X$.  It will
often be 
convenient to have a well-order of $X$ available;  we will then 
identify $X$ with the ordinal~$\lambda$.

We call a function $f: X \times X  \to  X  $ a
``pairing function'' if $f\on \{(x,y): x\not= y\} $
 is 1-1.    For the rest of the paper we
fix a pairing function~$pr$.   We will assume that the cardinality 
of the complement of the range of $pr$  is equal to the cardinality
of~$X$: $|X\setminus \ran(pr)|= \lambda $.

We fix a value~$0\in X$, and we will assume that $0$ is not in the
range of~$pr$.




When we consider terms in which several functions are nested, we may
write $f x $ or $g x y$ for $f(x)$ or~$g(x,y)$ to avoid too many 
parentheses. 

We identify $X^n$ with the set of functions from $\{1,\ldots, n\}$ to
$X$.   If $s\cap t = \emptyset$, $s \cup t = \{1,\ldots, n\}$, 
$a: s \to X$, $b: t \to X$, then $a \cup
b$ is in $X^n $.

If $\C \subseteq \O$ is a clone, we let $\oc n = \C \cap \oo n $.

\end{Notation}

\begin{Fact}\label{fact.zorn}
\begin{enumerate}
\item 
  If $f : X \times X  \to  X  $ is a
 pairing function, then there are unary functions~$g$, $g_1$, $g_2$
 such that  the function $(x,y)\mapsto g\circ f(g_1 x, g_2 y )$ is a
 bijection from $X\times X$ to~$X$. 
\item  If $\C \subseteq \O$, $\{pr\}\cup \oo 1  \subseteq 
 \C$, where $pr$ is any pairing function, 
then~$\C = \O$. [Use (1)]
\item 
If $\oo 1 \subseteq \C \subseteq \O$, $ pr \notin \C$, then the clones
which are maximal in 
$$ \{ \D: \ \C \subseteq \D \subseteq \O, pr \notin \D\}$$
are exactly the precomplete clones  extending~$\C$.   (Using Zorn's
lemma this easily implies that the interval 
 $[\oon  , \O]$  is dually atomic:  Every clone above~$\oo  1 $,
except for $\O$ itself, is contained in a precomplete one.)
\end{enumerate}
\end{Fact}

\begin{Remark} As we shall see in in section \ref{section.nobinary},
 we cannot relax the assumption
``$pr$ is 1-1 on $\{(x,y): x\not= y\}$'' in \ref{fact.zorn}(2) to 
``$pr$ is 1-1 on $\{(x,y): x< y\}$.'' 
\end{Remark}

\begin{Definition} \label{def.pol}
Let $I$ be any index set,   and $R \subseteq X ^I$.  Let~$f\in \oo n$. 

We say that $f $ respects $R$ iff:   
\begin{quote}
whenever $\bar \rho^1 = \langle \rho^1_i: i\in I\rangle,
 \ldots, \bar \rho^ n = \langle \rho^n_i: i\in I\rangle
$ are
all in~$R$, \\
then also $\langle f(\rho^1_i, \ldots, \rho^n_i) : i \in I \rangle\in
 R$.
\end{quote}

We let $\Pol R $ be the set of all functions respecting~$R$. 

\end{Definition}
We will usually be interested in the case where $R$ is a set of $n$-ary
functions on~$X$, i.e. $R \subseteq X^ {X^ n}$. 

\begin{Fact}\label{fact.pol}
The following observations follow easily from the definitions and from
the facts above. 
\begin{enumerate}
\item For any relation~$R$,  $\Pol R \subseteq \O$ is a clone. 
\item   If $\C$ is a clone, then $\C \subseteq \Pol \oc n$. 
\item If $\C$ is a clone and $\oc n \not= \oo n$, then 
$\Pol \oc n \subsetneq \O$.  
In fact, $(\Pol \oc n)^ {(n)} = \oc  n$. 
\item If $\C$ is a maximal clone and $\oc 1 \not= \oo 1$, then $\C =
\Pol \oc 1$. 
\item If $\C$ is a maximal clone and $\oc 1 = \oo 1$, then $\oc 2
\not= \oo 2 $, and  $\C = \Pol \oc 2$. 
\end{enumerate}
\end{Fact}

\begin{More Notation}  Let $\C$ be a clone on the set~$X$. 
We let $\tilde  \C$ be the set of all functions 
$\bar f: X^n\to X^k$ ($n,k>0$) such that   
  each function $\pi^k_i \circ \bar f$ is in~$\C$. 
\end{More Notation}

The ``closure under composition'' of the clone $ \C$ just means that 
$\tilde \C$ is closed under the usual notion of composition, i.e., 
whenever $\bar f: X^n\to X^m$ and  $\bar g: X^m\to X^k$ are in $\tilde \C$
 then also $\bar g \circ \bar f \in \tilde \C$.

\begin{Acknowledgment}
We are grateful to Lutz Heindorf for his thoughtful remarks on an
earlier version of the paper. 
\end{Acknowledgment}

\section{A new proof of Rosenberg's theorem}
\label{section.easy}

Let $X$ be an infinite set.     Rosenberg \cite{Ro76}
 has shown that there
are $2^{2^{|X|}}$ many precomplete clones on~$X$.  
Using transfinite induction 
he first constructs  $2^{2^{|X|}}$ many clones
with certain orthogonality properties and then shows that
they can be extended to pairwise different precomplete clones.

We give here an alternative proof of Rosenberg's theorem, utilizing
the well-known fact (see e.g.~\cite{CN2}) that on every infinite
set $X$ there are    $2^{2^{|X|}}$ ultrafilters.  
We will find an explicit  1-1 map from the ultrafilters to
precomplete clones.


\begin{Definition}
Let  $\I \subseteq \P(X)$ be a maximal  ideal.

  We define 
$$ \CI:= \bigcup_{n=1}^\infty \{ f\in \oo n: \ \forall A\in \I\,\,
f[A^n] \in \I\}$$
\end{Definition}

\begin{Fact}
\begin{enumerate}
\item $\CI \subsetneq  \O$
\item $\CI$ is a clone
\item If $f:X^k\to X$, $\ran(f)\in I$, then~$f\in \C_I$.  More
generally, if $A\in I$, $f:X^k \to X^n$
 and the range of $f$ is contained in~$A^n$,
then~$f\in \tilde \C$. 
\item $\I$ can be reconstructed from $\CI$ as 
$$ \I = \{A \subseteq X: 
\mbox{ For all $ f:X\to X$: If $\ran(f) \subseteq A$, then 
$ f\in \CI$}\},$$
so in particular the map $I \mapsto \CI$ is 1-1. 
\item 
   $\CI$ is a  precomplete clone, i.e.: For all $f\in \O \setminus \CI$ the
clone generated by $\CI\cup \{f\}$ contains all of~$\O$.
\end{enumerate}
\end{Fact}

\begin{proof}
Parts (1), (2) and  (3)  are clear.  We only check (4) and (5).

For (4), let 
$$ \I' = \{A \subseteq X: 
\mbox{ For all $ f:X\to X$: If $\ran(f) \subseteq A$, then 
$ f\in \CI$}\}.$$
By (3) above,  $\I \subseteq \I'$, so we check~$\I' \subseteq \I$. 
 Let~$A\notin \I$.  If $|A| \le |X\setminus A|$, then let $A_0:= A$,
otherwise we must have $|A| = |X|$, so we can write $A$ as a disjoint
union $A= A_0 \cup A_1$ with $|A_0|= |X| = |A_1|$,
$A_0 \notin \I$. \\
In either  case we have $A_0 \subseteq A$, $ |A_0| \le |X\setminus A_0|$, 
$A_0 \notin \I$. 
So there is a function 
$f:X \to X$ with $f[X]= f[X\setminus A_0] = A_0$, so $f\notin
\CI$ while $\ran (f)  \subseteq A_0 \subseteq A\in I$.   
Hence $A \notin I'$.

We now turn to the proof of (5). 
\begin{figure}
\bigskip

\setlength{\unitlength}{0.00053333in}
\begingroup\makeatletter\ifx\SetFigFont\undefined
\def\x#1#2#3#4#5#6#7\relax{\def\x{#1#2#3#4#5#6}}%
\expandafter\x\fmtname xxxxxx\relax \def\y{splain}%
\ifx\x\y   
\gdef\SetFigFont#1#2#3{%
  \ifnum #1<17\tiny\else \ifnum #1<20\small\else
  \ifnum #1<24\normalsize\else \ifnum #1<29\large\else
  \ifnum #1<34\Large\else \ifnum #1<41\LARGE\else
     \huge\fi\fi\fi\fi\fi\fi
  \csname #3\endcsname}%
\else
\gdef\SetFigFont#1#2#3{\begingroup
  \count@#1\relax \ifnum 25<\count@\count@25\fi
  \def\x{\endgroup\@setsize\SetFigFont{#2pt}}%
  \expandafter\x
    \csname \romannumeral\the\count@ pt\expandafter\endcsname
    \csname @\romannumeral\the\count@ pt\endcsname
  \csname #3\endcsname}%
\fi
\fi\endgroup
\newcommand{\dashlinestretch}{}
{  \renewcommand{\dashlinestretch}{30}
\begin{picture}(9395,6254)(0,-10)
\put(1137,3982){\ellipse{2250}{4500}}
\put(8262,3982){\ellipse{2250}{4500}}
\put(4887,2257){\ellipse{2250}{4500}}
\path(4137,607)(5637,607)
\path(7512,5632)(9012,5632)
\path(12,4282)(2262,4282)
\path(6012,607)(7137,2107)
\path(6312,232)(7437,1732)
\path(6012,907)(6012,607)(6237,607)
\path(7212,1732)(7437,1732)(7437,1507)
\path(2712,5632)(6612,5632)(6387,5407)
\path(6612,5632)(6387,5857)
\path(2112,2332)(3762,682)(3537,682)
\path(3762,682)(3762,832)
\put(4737,232){\makebox(0,0)[lb]{\smash{$A^k$}}}
\put(7962,5857){\makebox(0,0)[lb]{\smash{$B_1$}}}
\put(7887,4882){\makebox(0,0)[lb]{\smash{$B_0$}}}
\put(837,5032){\makebox(0,0)[lb]{\smash{$C_1$}}}
\put(1062,2857){\makebox(0,0)[lb]{\smash{$C_0$}}}
\put(12,5932){\makebox(0,0)[lb]{\smash{$X^n$}}}
\put(3837,4282){\makebox(0,0)[lb]{\smash{$X^k$}}}
\put(6987,832){\makebox(0,0)[lb]{\smash{$f$}}}
\put(6462,1507){\makebox(0,0)[lb]{\smash{$f^*$}}}
\put(4512,5782){\makebox(0,0)[lb]{\smash{$g$}}}
\put(2487,1432){\makebox(0,0)[lb]{\smash{$g_0'$}}}
\end{picture}
}

%
%
%
%

\bigskip
 \end{figure}

Call a function  $f$ ``conservative''  if it satisfies 
$f(a_1,\ldots, a_n)\in \{a_1,\ldots, a_n\}$ for all $a_1,\ldots,
a_n\in X$.  Clearly all conservative functions are in~$\CI$.

 Let~$f:X^k\to X$, $f\notin \CI$.   So there is some set $A\in \I$
with $f[A^k]\notin \I$.   Let $B_0 = f[A^k]$, $B_1 = X \setminus
B_0$.     So $B_0\notin \I$, $B_1\in \I$.  

Now let $g:X^n\to X$ be arbitrary.  We have to show that $g$ is in the
clone generated by $\CI $ and~$f$.   Pick two distinct elements $0,1$
in~$B_0$.      The function 
$$ H(x,y,z) = \KNUTHcases{y & if $x=0$\cr 
			z & if $x \not=0$\cr
		}$$
is conservative, hence in~$\CI$. 

Let  $C_0 = g^{-1}[B_0]$,  $C_1 = g^{-1}[B_1]$, and define two
``approximations''~$g_0$, $g_1$ to $g$ as follows: 
$$ g_0(\bar x) = 
\KNUTHcases{ g(\bar x) & if $\bar x \in C_0$\cr
            0           & if $\bar x \in C_1$\cr
          }
\qquad \qquad 
 g_1(\bar x) = 
\KNUTHcases{     0           & if $\bar x \in C_0$\cr
	g(\bar x) & if $\bar x \in C_1$\cr
          }
$$
Let $\chi(\bar x) = 0$ if $\bar x\in C_0$, $\chi(\bar x) = 1$
 if $\bar x\in C_1$.   By definition of~$H$, 
$g(\bar x) = H(\chi \bar  x, g_0 \bar x,
 g_1 \bar x)$, so all we have to show is that $H, \chi, g_0, g_1$ are all
 in the clone generated by $\CI$ and~$f$.    We already know that 
\begin{enumerate}
\item 
$H \in \CI $  (because $H$ is conservative), 
\item $ g_1\in \CI$ (because the range of $g_1$ is in~$\I$),
\item $ \chi\in \CI$ (because $ \chi $ takes only 2 values)
\end{enumerate}
It remains to show $ g_0\in cl(\CI \cup \{f\})$. 

Let $f^*: B_0\to A^ k$ be an ``inverse'' of~$f$, i.e., 
$$ \forall b \in B_0 : \ f( f^ *(b)) = b$$
($f^ * \on B_1$ can be   arbitrary  function with range~$ \subseteq A^
k$.)\\
Define $g_0': X^ n \to X^ k$ by $ g_0' (c) = f^ *(g_0(c))$.  Note that
the range of $g_0'$ is $ \subseteq A^ k$, $ A \in I$, so~$ g_0' \in
\tilde \CI$.     

Now we have, for all $\bar c\in X^n$, 
 $ g_0(\bar c) = f(f^*(g_0(\bar c)) = f(g_0'(\bar c))$, so $ g_0\in
cl(\CI\cup \{f\})$. 

\end{proof}

%
%
%
%
%
%

\begin{Conclusion}  On any infinite set $X$ there are exactly
$2^{2^{|X|}}$ many precomplete clones. 
\end{Conclusion}
\begin{proof}  The upper bound follows from $|\O| = 2^{|X|}$. For the
lower bound: it is known that there are
$2^{2^{|X|}}$ many maximal ideals, 
and we have just shown that the function $I\mapsto \CI$ maps them 
injectively 
to precomplete clones. 
\end{proof}

\section{Almost unary clones}
\label{section.nobinary}

In this section we will consider clones on an infinite set  $X$ of
regular cardinality.  We will call a set ``small'' if its cardinality
is smaller than the cardinality of~$X$, and we will say that there are
``few'' objects with some property if the set of those objects is
small. 

For example, $X$ is countable, then ``small'' will mean ``finite''.
If $X$ has cardinality ${\aleph_1} $, then  ``small'' will mean
``finite or countably infinite''.

[With this notation, the property ``$X$ has regular cardinality'' can
be rephrased as ``$X$  cannot be written as a union of few
 small sets'']

\begin{Definition}\label{def.almost}
Let~$g: X^n \to X$.     We say that $g$ is {\em
almost unary} iff there is a function $G$ which is defined on~$X$,
each $G(x) $ a small subset of~$X$, such that for some~$k$: 
 $$ \forall (x_1,\ldots, x_n)\in X ^n : g(
x_1,\ldots, x_n) \in G(x_k)$$

If $X$ itself is a cardinal, then we can equivalently say: $g$ is
almost unary iff:  for some~$k$, $G: X\to X$, 
 for all $x_1,\ldots, x_n\in X$:  $ g(x_1, \ldots, g_n)\le G(x_k)$. 

\end{Definition}
\begin{Definition}
Let $\U  \subseteq \O$  be the set of all almost unary functions.
\end{Definition}

In definition \ref{def.heavily} we will call functions in $\oo 2 
\setminus \U$ ``heavily binary''.

\begin{Definition}
Let $\hat \U:= \Pol \H $  (see \ref{def.pol}). 

That is, a function $f\in \oo n$ is in $\hat \U$ iff 
$$ \forall g_1,\ldots, g_n\in \H: f(g_1,\ldots, g_n)\in \H$$
where $f(g_1,\ldots, g_n)$ is the function 
$(x,y)\mapsto f(g_1(x,y),\ldots, g_n(x,y))$. 

Note that $\hat   \U \cap \oo 2  = \U\cap \oo 2 $, and~$\U \subseteq
\hat \U$.
\end{Definition}

\begin{Example} \label{example.med}
Let $X= \lambda $ be a cardinal, so 
the small subsets of $X$ are exactly the bounded subset of~$\lambda $.
\begin{enumerate}
\item The function $\min$ is almost unary: $\min \in \H$
\item the function  $\max$ is not almost unary.
\item The median function~$med$, defined by 
$$ med(x,y,z) = \max(\min(x,y),\min(y,z),\min(x,z)) 
$$
is not almost unary,  but it is easy to check that 
 $med$ respects all almost unary functions, so $med 
 \in \hat \U \setminus \U $. 
\item Let $pr_\Delta $ be defined by 
$$
pr_\Delta (x,y) = 
\KNUTHcases{  pr(x,y)          & if $x>y$\cr
0	& otherwise
          }
$$
(where $pr$ is a pairing function). Then $pr_\Delta\in \U$. 
\end{enumerate}
\end{Example}

The following was already observed by Davies and Rosenberg
\cite{DR85}. 
\begin{Conclusion}\label{dr}
Assume $\C\in [\O^{\langle 1\rangle},\O]$.   If $pr_\Delta \in \C$
(see~\ref{example.med}), and if $\C$ contains a binary function not in
$ \H$, then~$\C = \O$.

Hence, $\Pol \H $ is an example of a  precomplete clone containing all
unary functions. 
\end{Conclusion}

\begin{proof}  Let $p_1$, $p_2$: $X\to X$ be two 1-1 functions such
that the ranges of $p_1$, $p_2$, $pr$ are disjoint. 
  Since $\C$ contains a function which is not almost
unary, there is some  $H\in \oc 2 $ with 
 $H(x,p_2 0)=x  = H(p_1 0,x)$ for all $x$ in the range of~$pr$. 
  Then the function  $$(x,y)\mapsto
H(\, p_1( pr_\Delta (x,y)) \, ,\,  p_2 ( pr_\Delta (y,x))\, )$$
is a pairing function. 
\end{proof}

(We will meet a similar argument again in the proof of
\ref{lemma.HH}.) 

We now show a kind of  converse to this theorem: $\Pol\H$ is the {\em
unique}
precomplete clones which which contains all unary
functions and only ``almost unary'' binary functions. 

\begin{Theorem}\label{thm.u2}
Assume that $\C \subseteq \O$ is a precomplete clone, $\oo 1 \subseteq
\C$, $\oc 2  \subseteq \H$. \\
Then $\C = \Pol \H $. 
\end{Theorem}

We will prove this theorem below.  We start by investigating which
coordinates are responsible for a function having a large  range.

\begin{Definition}  Let $g\in \oo n$.  We define a set $S_g $ of
subsets of $\{1,\ldots, n\}$ as follows. 
$$ S_g = \{s \subseteq \{1,\ldots, n\}: 
\exists \bar a \in X ^{\{1,\ldots, n\}\setminus s}: \ 
\bigl|\{g(\bar a \cup \bar x): \bar x\in X ^s\}\bigr| =
\bigl |X  \bigr|\, \}
$$

(Here we write   $X^s $  for the set of all functions from $s$
 to~$X$.) 
\end{Definition}

\begin{Lemma}\label{lemma.st}
Assume $cl(g)^{(2)} \subseteq \U$. Then 
$$ \forall r,t\in S_g:  r \cap t \not= \emptyset$$
\end{Lemma}

\begin{proof}
   Choose $r$ and $t$ in~$S_g$ with $r\cap t  = \emptyset $. 
 Using unary
functions, we will construct a binary  function in $cl(g)$ which is not
in~$\H$. 

Let $s:= \{1,\ldots, n\} \setminus (r\cup t)$, so 
$ \{1,\ldots, n\} = r \dot{\cup} s\dot{\cup} t$. 
 So there is some $\bar a\in
X^{s\cup t}$ and a sequence $(\bar x^\alpha: \alpha \in X)$ of elements of
$X ^r$ such that all values $g(\bar a\cup \bar x^\alpha)$ are different.
Similarly, 
 there is some $\bar b\in
X^{r\cup s}$ and a sequence $(\bar y^\beta: \beta \in X)$ of elements of
$X ^t$ such that all values $g(\bar b\cup \bar y^\beta)$ are different.

Now for $\ell=1,\ldots, n$ define functions $h_\ell$ as follows: 
Fix some
element~$0\in X$. 
$$
h_\ell(\alpha,\beta ) = \KNUTHcases {
                \bar x_\alpha(\ell) & if $\ell\in r$, $\alpha \not= 0$\cr
                \bar b(\ell) & if $\ell\in r$, $\alpha = 0$\cr
                \bar a(\ell) & if $\ell\in s$, $\alpha \not= 0$\cr
                \bar b(\ell) & if $\ell\in s$, $\alpha = 0$ \cr
                \bar y_\beta(\ell) & if $\ell\in t$, $\beta\not=0$\cr
                \bar a(\ell) & if $\ell\in t$, $\beta = 0$\cr
}
$$
Formally, the functions $h_\ell$ are in~$\oo 2$, but each of them
is essentially unary: $h_\ell(\alpha, \beta)$  depends only  on
$\alpha$ for $\ell \in r\cup s$, and only on $\beta $ for $\ell\in
t$.  This implies that $h_\ell\in \O^{\langle 1\rangle}$.

Now the function $F = g(h_1, \ldots, h_n)$, i.e., 
$F (\alpha,\beta) = 
g(h_1(\alpha,\beta ), \ldots,  h_n(\alpha, \beta))$, will be in $\oc 2
$ but not in~$\H$, since the values $F(\alpha,0) = g(\bar x_\alpha
\cup \bar a)$ 
are all different, as are the values $F(0, \beta ) = g(\bar y_\beta
\cup \bar b)$.  
\end{proof}

The previous lemma will allow us to relate any ``almost unary'' clone
to $\Pol\H$: 

\begin{Lemma} \label{lemma.d}
Assume $\oo 1 \subseteq \C$, $\oc 2 \subseteq \H$.   Then 
$\C \subseteq \Pol \H$.    That is: whenever $d_1,\ldots, d_n\in \H$, 
$g\in \oc n$, then also $f:= g(d_1,\ldots, d_n)\in \H$. 
\end{Lemma}
\begin{proof}
Since each $d_\ell\in \H$, we can find a decomposition $\{1,\ldots,
n\} = r \cup t$, $r\cap t= \emptyset$, and a function $D$ mapping each
$\alpha \in X$ to a small subset $D(\alpha ) \subseteq X$ such that:
\begin{enumerate}
\item For all $\ell\in r$, all $\alpha,\beta\in X$:
$d_\ell(\alpha,\beta )\in
D(\alpha)$.
\item For all $\ell\in t$, all $\alpha,\beta\in X$:
$d_\ell(\alpha,\beta)\in D(\beta )$.
\end{enumerate}
By the previous lemma, we cannot have both $r$ and $t$
in~$S_g$, so wlog assume~$t\notin S_g$. 

Now fix any element~$0\in X$.  We will show that the set
 $\{f(0,\beta): \beta \in X\}$ is small.  

Consider  $f(0,\beta) = g(d_1(0,\beta), \ldots, d_n(0,\beta))$.
Identifying $X^n$ with $X^{\{1,\ldots, n\}}$, we can write the tuple 
$(d_1(0,\beta), \ldots, d_n(0,\beta))$ as $a_\beta  \cup y_\beta$,
$a_\beta \in X^r$, $y_\beta \in X^t$. Now note that
 for
 $\ell\in r$ we have $d_\ell(0, \beta) \in D(0)$,
 so $a_\beta \in D(0)^r$, which is a small set. 

Hence 
$$ \{f(0,\beta): \beta \in X\} \subseteq \{g(a\cup y): a\in D(0)^s, y\in
X^t\}$$
For each fixed $a\in D(0)^s $ the set  
$\{g(a\cup y): y\in X^t\}$
is small (since $t\notin S_g$), so,
since $D(0)$ is small, 
 also 
$$
\{g(a\cup y): a\in D(0)^s, y\in
X^t\} = \bigcup_{a\in D(0)^r} \{ g(a\cup y): y\in X^t\}$$
is small. 
\end{proof}

\subsection*{Proof of theorem \ref{thm.u2}}

Assume $\oo 1 \subseteq \C$, $\oc 2 \subseteq \H$, and assume that $\C$
is precomplete.     Then by lemma \ref{lemma.d}, we have 
$\C \subseteq \Pol \H$.  But since $\C$ is maximal, we must have $\C
  = \Pol \H$.

%
%
%

\section{Successors of regulars}
\label{section.many}

We fix a set $X$ of regular cardinality~$\lambda $, and  for simplicity we
write~$X = \lambda $. We fix a pairing
function $pr:\lambda \times \lambda \to \lambda $ as in \ref{notation}. 

We will use the following combinatorial principle $Pr(\lambda,\mu)$:  

\begin{quote}
There is a symmetric function $c:\lambda \times \lambda 
 \to \mu$ with the following
anti-Ramsey property: 
\\
For all sequences $(a_i:i< \lambda )$ of
pairwise disjoint finite subsets of $\lambda $, and
for all $c_0\in \mu$:
$$ \mbox{ there are $i<j<\lambda $ such that  
$c \on (a_i\times a_j) $ is constant with value $c_0$}$$
\end{quote}
(See section \ref{section.app} for background)

We fix a function $c$ witnessing the above statement. 


\begin{Definition}\label{def.fa}
For any  $ A \subseteq \mu$ we define a function $F_A  :\lambda \times
\lambda \to \lambda  $ as follows: 
$$ F_A(\alpha, \beta) = 
\KNUTHcases{ \max(\alpha, {\beta})       & if $\alpha=0$ or $\beta = 0 $
					or $\alpha = \beta $\cr
	pr(\alpha,\beta )	& if $c(\alpha,\beta)\in A$ \cr
	0			& otherwise\cr
}$$
\end{Definition}

\begin{Fact}\label{fact.many}
If $A \cup B = \mu$, then $cl(F_A, F_B) = \O$. 
\end{Fact}

\begin{proof}
We will show how to construct a pairing function from 
$F_A$ and~$F_B$. 

Define $$ pr'(\alpha,\beta) = F_A(F_A(\alpha,\beta),
F_B(\alpha,\beta))$$

We claim that for all distinct $\alpha, \beta > 0$: $pr'(\alpha,\beta)
= pr(\alpha, \beta)$. \\
Indeed, if  $c(\alpha,\beta)\in A\cap B$, then 
$$
 pr'(\alpha,\beta) = F_A(pr(\alpha,\beta), pr(\alpha,\beta) )
 = pr(\alpha,\beta),$$
if   $c(\alpha,\beta)\in A \setminus B$, then 
$$
 pr'(\alpha,\beta) = F_A(pr(\alpha,\beta), 0) = pr(\alpha,\beta)$$
and if   $c(\alpha,\beta)\in B \setminus A $, then 
$$
 pr'(\alpha,\beta) = F_A(0, pr(\alpha,\beta)) = pr(\alpha,\beta).$$

Hence $pr'$ is a pairing function. 
\end{proof}

\begin{Main Lemma}\label{lemma.many}
Assume that $A \not \subseteq  B_1 \cup \cdots B_k$.  Then 
$$ F_A \notin cl(F_{B_1} , \ldots, F_{B_1} ).$$
\end{Main Lemma}
We will prove this lemma below, but first we will show how it can be
used. 

\begin{Definition}
We say that ${\mathscr A} = (A_i:i\in I)$ is an
 {\em independent} family of subsets of~$X$, if 
 every
nontrivial Boolean combination of sets from ${\mathscr A}$ is nonempty, 
i.e.:   
\\ 
Whenever $J_0$ and $J_1$ are finite disjoint subsets of~$I$, then 
$$ \bigcap_{i\in J_0} A_i \ \cap \bigcap_{i\in J_1} (X\setminus A_i) \
\ \not= \ \ \emptyset$$
\end{Definition}
The following theorem of Hausdorff is well known: 
\begin{Theorem}
If $|X| = \mu$, then there is an independent family 
 ${\mathscr A}= 
(A_i:i \in I)$ of subsets of $X$ with $|I| = 2 ^ \mu$. 
\end{Theorem}

\begin{proof} See \cite[Chapter VIII, exercise A6]{Ku83} or 
\cite[Example 9.21]{Ko89}. 
\end{proof}

\begin{Theorem}\label{theorem.many}
 Assume $Pr(\lambda , \mu)$.  Then there are at least 
$ 2^{ 2^  \mu} $ many precomplete clones above the unary functions on the
set~$\lambda  $.  (Hence:  If $\lambda = \kappa ^+ $, $\kappa $
regular, then there are  $2^{2^\lambda}$ many precomplete clones above
$\oo 1$.)
\end{Theorem}

\begin{proof}
Let $(A_i: i \in 2^ \mu)$ be an independent family of subsets of~$\mu$.
  Write $-A_i$ for
$ \mu  \setminus A_i$.   For each $J
\subseteq 2 ^ \mu$ we let 
\begin{quote}
$ \C_J = $ the clone generated by 
$ \{ F_{A_i}:i \in J\} \cup \{ F_{-A_i}:i \notin J\}\cup \oo 1  $. 
\end{quote}

We will now show that 
\begin{enumerate}
\item $ \C_J \not= \O$, for all $ J \subseteq 2^ \mu $
\item Whenever $J_1 \not= J_2$, then $ \C_{J_1} \cup 
 \C_{J_2}$ already generates~$ \O$. 
\end{enumerate}
This will conclude the proof, because (1) together with fact
\ref{fact.many} implies that 
each $\C_J$ can be extended to a precomplete clone, and (2) implies that
no single precomplete clone can contain  $ \C_{J_1} \cup 
 \C_{J_2}$ for distinct~$J_1$,~$J_2$. 

Proof of (1):   Wlog there is some~$i\notin J$.  By independence, 
 $A_i$ cannot be
covered by any finite union from $\{ A_j: j \in J\} \cup  
 \{ -A_j: j \notin J\}$. So by the lemma, $ F_{A_i}$ is not in the
clone~$\C_J$. 

Proof of (2):  If $J_1\not= J_2$, then there is wlog some $i\in J_1
\setminus J_2$.  Now $F_{A_i}\in \C_{J_1}$,  $F_{-A_i}\in \C_{J_2}$,
and 
by fact~\ref{fact.many}, $\{F_A, F_{-A}\}$ generates~$\O$. 
\end{proof}

We now prepare for the proof of the main lemma \ref{lemma.many}. 
Our situation is the following:  We have a function $c$ witnessing
$Pr(\lambda,\mu)$.  Using $ c$ and our fixed pairing
function $pr$ we have defined functions
$F_A: \lambda \times \lambda \to \lambda $ for every $A \subseteq
\mu$ in \ref{def.fa}. 
   We are given sets   $A, B_1,\ldots, B_k \subseteq \mu$, $A
\not \subseteq 
B_1\cup \cdots \cup B_k$.  Pick $c_0\in A\setminus (B_1\cup \cdots
\cup B_k)$. 

     We want to show that $F_A\notin
cl(F_{B_1},\ldots, F_{B_k})$, i.e., the functions 
$ F_{B_1},\ldots, F_{B_k}$, together with all unary functions, do not
generate~$F_A$.

\begin{Definition}
  ``Terms'' over $\lambda  $ are defined inductively as
follows: 
\begin{enumerate}
\item The formal variables~$\x$, $\y$ are terms, as well as
every element of~$ \lambda $. 
\item If $ \sigma $ is a term, $f: \lambda  \to  \lambda $ a unary function, then $
  (f,\sigma)$ is a term. 
\item If $ \sigma_1 $ and $\sigma_2 $ are terms, $1 \le i\le k$, 
 then $(F_{B_i},\sigma_1,\sigma_2)$ is a
term. 
\end{enumerate}
Every term  $\tau  $ induces (in the obvious way) a function
$\tau: \lambda \times  \lambda  \to  \lambda $ which is in $cl(\{F_{B_1},\ldots  , F_{B_k}\})$. 
Conversely, every function in $cl(\{F_{B_1},\ldots  , F_{B_k}\})$ is
represented by a term.

We call a term ``constant'' if it is an element of~$ \lambda $, and we call a
term $\x$-unary  if $\y$  does not appear in it, similarly for
$\x$-unary.     A term is unary if it is $\x$-unary or $\y$-unary.
(By definition, the constant terms are both $\x$-unary and
$\y$-unary.)

\end{Definition}

For the following discussion, fix a term~$\tau_0$.   Our aim is to
find a large set on which all subterms of $\tau_0$ behave like unary
functions.    We will first explain how to find (terms for) 
these unary functions, and then we show 
they are indeed realized on some large
set.

\begin{Definition}  Let~$S \subseteq  \lambda $.   For any term $\tau$
 we will try to define a unary term~$\tau^S$.
Whenever $\sigma^S$ is undefined for a subterm $\sigma $ of~$\tau$,
 then also
$\tau^S$ will be undefined.   Our definition proceeds by induction on
the structure of~$ \tau$.  ``$B$'' will stand for any of the sets
$B_1$, \dots,~$B_n$. 
\begin{enumerate}
\item   $\tau = \x$ or $\tau=\y$ or $\tau=c \in  \lambda $. \\ In this case,
$\tau^S = \tau$. 
\item  $\tau = (f, \sigma)$, and $\sigma^S = c\in  \lambda $. \\ In this case,
$\tau^S$ is also a constant, namely:~$f(c)$.  
\item  $\tau = (f, \sigma)$,  $\sigma^S = (g, \x)$. 
\\ If $f\circ g$ is 1-1 on~$S$, then $\tau^S:= (f\circ g, \x)$. 
\\ If $f\circ g$ is constant with value $d$ on~$S$, then~$\tau^S:= d$.
\\ If $f\circ g$ is neither 1-1 nor constant,
 then $\tau^S $ will be undefined. 
\item $\tau = (F_B, \sigma_1, \sigma_2)$, and $\sigma_1^S$ and
$\sigma_2^S$ are constant (say, with values $c_1$ and $c_2$): 
\\ In this case we let $\tau^S:= F_B(c_1, c_2)$. 

\item $\tau = (F_B, \sigma_1, \sigma_2)$, $\sigma_1^S = (f, \x)$,
$\sigma_2 = d$ (a constant). \\
If the function $h: x \mapsto  F_B(f(x), d)$ is 1-1 or constant (say, with
value $=c $) on~$S$,
 then we let 
$ \tau^S:= (h, \x)$ or~$=c$, respectively.   (If $h$ is neither constant
nor 1-1 on~$S$, then  $\tau^S$ is again undefined.) 
\item $\tau = (F_B, \sigma_1, \sigma_2)$, and $\sigma_1^S= ( f_1,\x)$ 
$\sigma_2^S= ( f_2,\x)$. 
\\If the function $h:x \mapsto F_B(f_1(x), f_2(x))$ is 1-1 or constant (say,
with value $=d$),
then we let $\tau^S = (h,\x) $ or~$ d$, respectively. (Otherwise,
$\tau^S$ is again undefined.) 
\item $\tau = (F_B, \sigma_1, \sigma_2)$, and $\sigma_1^S = (f_1,\x)$, 
$\sigma_2^S =(f_2,\y)$. \\
We let $\tau^S := 0 $.
{\bf This is the crucial  case  of our
definition.} 
\item Repeat all the above items with  $\x$ and~$\y$ interchanged,
 and/or $\sigma_1$ and $\sigma_2$ interchanged. 
\end{enumerate}
\end{Definition}

\begin{Fact} Whenever $\tau^S$ is defined, then $\tau^S$ is either
constant, or of the form $(f,\x)$ or~$(f,\y)$, where $f$ is 1-1 on
$S$. 
\end{Fact}

\begin{Fact} \
\begin{enumerate} 
\item
If  $\tau^S$ is defined and $S' \subseteq S$ , then $\tau^{S'}$ is
defined. 
\item 
Fix a finite set $T$ of terms which is closed under
subterms. Then:  for every  set $S$  of regular infinite cardinality
there a set $S' \subseteq S$ of the same cardinality such that:
\begin{quote}
For all~$\tau \in T$, $ \tau^{S'}$ is well-defined.
\end{quote}


\end{enumerate}
\end{Fact}
\begin{proof} Proceed by  induction on the complexity of the terms.
We have to thin out the set $S$ finitely many times in order to make
finitely many functions 1-1 or constant. 
\end{proof}

\begin{Lemma}\label{lemma.submain}
Assume that $\tau^S$ is defined,~$|S| = \lambda $.
  Then there are $\alpha < \beta $ 
in $S$ such that  $\tau(\alpha,\beta) = \tau^S(\alpha,\beta)$ and 
$c(\alpha, {\beta} ) = c_0$. 
\end{Lemma}
\begin{proof}
Let $T$ be the set of subterms of $\tau$ (including $\tau$ itself). 
Collect all the 1-1 functions appearing in $\sigma^S$ for $\sigma \in
T$, i.e.: 
$$ \F:= \{ f:\  \exists \sigma\in T\  \sigma^S=(f,\x) 
\mbox{ or } \sigma^S=(f,\y) \}$$
The set  $\F$ is finite, the identity function is in~$\F$, 
 and all functions in $\F$ are 1-1. We may
thin out the set $S$ so that the family 
$$ (\{f(\alpha): f\in \F\} : \alpha \in S)$$
is pairwise disjoint.  So since $c$ witnesses $Pr(\lambda,\mu)$, 
 we can find $ \alpha < \beta $ such that

For all $f,g\in \F$:   $c(f(\alpha), g(\beta)) = c_0$ (and
$f(\alpha)\not=g(\beta)$). 

This implies $F_{B_i}(f(\alpha), g(\beta))  = 0$. 

Now we can prove by induction on the complexity of the subterms $
\sigma $ of $\tau$ that $\sigma^S(\alpha, \beta) =
\sigma(\alpha,\beta)$. 
\end{proof}

\subsection*{Proof of lemma \ref{lemma.many}}
Let $c_0 \in A  \setminus (  B_1 \cup \cdots B_k)$, and let $\tau$ be
a term.   We will find $\alpha,\beta$ such that
$\tau(\alpha,\beta)\not= F_A(\alpha,\beta)$. 

We can find a set $S$ such that  $\tau^ S$ is defined.  Let $\F$ be
again the finite set of 1-1 functions used in defining~$\tau^S$.
  We can thin out the set $S$ such
that for all $f\in \F$: 
$$ \forall \alpha, \beta \in S:  
\ \alpha \not= \beta \ \Rightarrow \ 
f(\alpha) \not =  pr(\alpha, \beta) \not= f(\beta)$$

[Why?  For each such $f\in \F$ define a partial  function $\bar f$ such that
$\bar f(\alpha) = \beta $ whenever $f(\alpha) = pr(\alpha, \beta )$, $
\alpha \not= \beta $.
$\bar f $ is well-defined, since $pr$ is a pairing function.   We can
thin out $S$ to get:
 $ \forall \alpha\in S: \bar f(\alpha)\notin S$. This is sufficient.]

Now thin out $S$ such that  $ \forall \alpha \in S$: $f(\alpha) \notin
S$ or $f(\alpha)=\alpha$, and that none of the finitely many
constants appearing as $\tau^S$  is
equal to $pr(\alpha,\beta)$ for $\alpha,\beta\in S$.

By lemma \ref{lemma.submain},
 we can find $ \alpha < \beta $ with $\tau(\alpha,\beta)
= \tau^S(\alpha,\beta)$, and 
$c(\alpha, \beta) = c_0$.  Now we have $F_A(\alpha,\beta)  =
pr(\alpha,\beta)$ (as $c(\alpha,\beta) = c_0\in A$).  On the other
hand, $\tau^S$ is either constant or of the form $(f,\x)$ or $(f,\y)$
for some~$f\in\F$. So $\tau^S(\alpha ,\beta ) \not=
F_A(\alpha,\beta)$. 

This concludes the proof of lemma \ref{lemma.many} and hence also of
theorem \ref{theorem.many}.

\section{Weakly compact cardinals} \label{section.wc}

In this section we deal with clones on infinite sets whose  cardinality
$\lambda $ satisfies $ \lambda \to (\lambda)^2_2$ (so either $\lambda
= {\aleph_0}$ or $\lambda$ is weakly compact).

Recall that $\lambda \to (\lambda)^2_2$ implies 
$$ 
\lambda \to (\lambda)^n_k$$
for all $n,k< \omega $, i.e.:   Whenever $h:[\lambda] ^n \to \{1,
\ldots, k\}$, then there is a subset $S \subseteq \lambda $, $|S|=
\lambda $  such that
$h\on [S]^n$  is constant.

\begin{Definition}\label{def.heavily}
Let $H: \lambda ^ n \to  \lambda$.  
\begin{enumerate}
\item 
  We say that ``$H$ depends on the
$k$-th coordinate'' iff  there is
 $ (a_1, \ldots, a_{k-1}, a_{k+1},\ldots,  a_n)$ such that  the set 
$$ \{ H(a_1, \ldots, a_{k-1}, x,\ldots,  a_n):  x \in \lambda \} $$
has more than one element. 
In this case we may also write $H$ symbolically as $H(\x_1,\ldots,
\x_n)$ and  say ``$H $ depends on $\x_k$''.   For $n=2$ we may also
say ``$H(\x,\y)$  depends on $\x$''  or ``\dots\ on $\y$''. 
\item 
We say that $H(\x_1,\ldots, \x_n)$ {\em depends heavily} on the $k$-th
  coordinate (or: ``on $\x_k$'') iff 
  there is an $n-1$-tuple 
 $ (a_1, \ldots, a_{k-1}, a_{k+1},\ldots,  a_n)$ such that  the set 
$$ \{ H(a_1, \ldots, a_{k-1}, x,\ldots,  a_n):  x \in \lambda \} $$
has $\lambda $ many  elements. 
\item We say that $\C \subseteq \O$ is ``heavily binary'' if there
exists $H(\x,\y)\in \C$, which  depends
heavily on $\x$ and which also depends  heavily on~$\y$.  
\end{enumerate}
Thus, the  functions which are not ``heavily binary'' are exactly the 
``almost unary'' functions of definition~\ref{def.almost},
 and the heavily binary clones   are exactly those $\C \subseteq \O $
which satisfy $\C^{(2)} \not\subseteq \U$. 
\end{Definition}

\begin{Example}\label{example.h}
Let $H:\lambda \times \lambda \to \lambda $ be a function satisfying 
$$(*) \qquad \qquad \qquad 
 \forall \alpha > 0:  \ H(\alpha, 1) = \alpha  = H(0,\alpha)$$
[E.g., the $\max$ function has this property.] \\
Then $H$ depends heavily on $\x$ and~$\y$. 

Conversely, if $\C$ is a clone containing all unary functions and at
least one heavily binary function, then $\C$ contains a function $H$
satisfying $(*)$ above. 
\end{Example}

The following example shows that there are nontrivial heavily binary
clones above~$\oo 1 $. 
\begin{Example}\label{example.pq}
We will write $[X]^{<n}$ for the family of 
 subsets of $X$ of size~$<n$, and we
will write  $[X]^{<{\aleph_0} }$ for the family of finite  subsets of
~$X$. 
\begin{enumerate}
\item 
Let $\Q$ be the set of all functions $f\in \O$ such that: 
\begin{quote}
for some~$n$, $f: X^n \to X$, \\ 
and there is a function $Q:X \to [X]^{<{\aleph_0} }$, 
$$ \forall x_1\cdots x_n: \ f(x_1,\ldots, x_n) \in Q(x_1)\cup \cdots
\cup Q(x_n)$$
\end{quote}
\item 
Let $\P$ be the set of all functions $f\in \O$ such that: 
\begin{quote}
for some~$n$, $f: X^n \to X$, \\ 
and there is some $k$ and  a function $P:X \to [X]^{<k }$, 
$$ \forall x_1\cdots x_n: \ f(x_1,\ldots, x_n) \in P(x_1)\cup \cdots
\cup P(x_n)$$
\end{quote}
\end{enumerate}
Then:
\def\itm #1 {\item[{(#1)}]}
\begin{enumerate}
\itm A  $\P$ and $\Q$ are clones.
\itm B  $\P \subseteq \Q \subseteq \O$. 
\itm C $\P$ contains a heavily binary function, as well as all unary
functions.  
\itm D  If $X$
is finite, then trivially~$\P=\Q=\O$. 
\itm E If $X$ is countably infinite, then $\P \subsetneq \Q = \O$. 
\itm F If $X$ is uncountable, then $\P \subsetneq \Q \subsetneq \O$. 
\end{enumerate}

We leave the verification of this fact to the reader. 
\end{Example}

\begin{Theorem} Assume that $ \lambda \to (\lambda) ^ 2 _ 2$, i.e., $
\lambda $ is weakly compact or $ \lambda = {\aleph_0}$.   Then there
is a unique precomplete clone which contains all unary functions and 
 is heavily binary. 

By example \ref{example.pq} there are nontrivial heavily binary clones
above $\oo 1$, so by fact~\ref{fact.zorn}(3) there must be at least
one precomplete such clone. So  it is enough to show the following: 
 Whenever~$ \C_1$, $\C_2$ are heavily binary 
clones on $ \lambda $, $
\oo 1  \subseteq \C_1 \cap \C_2$, then $cl( \C_1 \cup \C_2) = \O$ implies
$\C_1 = \O$ or $\C_2 = \O$. 
\end{Theorem}

To make the proof clearer, we need a few definitions and lemmas.

\begin{Definition}
For $S \subseteq \lambda$, let 
$$ \Delta_S = \{(\alpha,\beta)\in S\times S:  \alpha > \beta \}
\qquad     \qquad     
 \nabla_S = \{(\alpha,\beta)\in S\times S:  \alpha < \beta \}
$$
We let $\nd_S:= \nabla_S \cup \Delta_S = \{ (\alpha,\beta)\in S\times
S: \alpha \not= \beta \}$. 
\end{Definition}

\begin{Definition} 
For $\bar \alpha =(  \alpha_1,\alpha_2,\alpha_3, \alpha_4)\in
\lambda^4$,
$\bar \beta =(  \beta_1,\beta_2,\beta_3, \beta_4)\in
\lambda^4$ we define $ \bar \alpha \sim \bar \beta $ iff $ \forall
i,j\in \{1,2,3,4\} : (\alpha _i< \alpha_j \Leftrightarrow
\beta_i<\beta_j)$. 
\end{Definition}

\begin{Definition}\label{def.canonical}
Let $ F: \nd_S\to \lambda $. We say that $ F$ is {\em canonical}
 on $S$ iff: 
\\
For all $ \bar \alpha \sim \bar \beta$: If
$F(\alpha_1,\alpha_2)<F(\alpha_3,\alpha_4)$, then 
$F(\beta_1,\beta_2)<F(\beta_3,\beta_4)$.

[This also implies: 
For all $ \bar \alpha \sim \bar \beta$: If
$F(\alpha_1,\alpha_2)=F(\alpha_3,\alpha_4)$, then 
$F(\beta_1,\beta_2)=F(\beta_3,\beta_4)$.]

\end{Definition}
\begin{Fact}\label{fact.reader}
 If $ \lambda \to (\lambda)^2_2$, then 
$$ \forall S \in [\lambda ] ^ \lambda \ \exists S' \in [S]^\lambda: \ 
\mbox{$F$ is canonical on $ S'$}$$
\end{Fact}

The proof uses the partition relation $\lambda \to (\lambda )^4_n$ for
some large~$n$.   We leave the details to the reader. See also 
fact~\ref{fact.allfinite}. 


\begin{Lemma}\label{lemma.canonical}
Assume that $ F$ is canonical on~$\nd$. Then: 
\begin{enumerate}
\item $F\on {\Delta}$  satisfies one of the
following properties:  
\begin{itemize}
\item $F\on \Delta $ is 1-1 [typical examples: $pr$, $pr_{\Delta}$.]
\item $F\on {\Delta} $ depends injectively 
 on the first coordinate:  $F(x,y)=
g(x)$ for some 1-1 function~$g$. [typical examples: $\pi^2_1$, $\max$] 
\item $F\on {\Delta} $ depends injectively on 
 the second  coordinate:  $F(x,y)=
g(y)$ for some 1-1 function~$g$. [typical examples: $\pi^2_2$, $\min$] 
\item $F\on {\Delta} $ is constant. 
\end{itemize}
\item Similarly for~$F\on \nabla $. 
\item If at least one of $F\on \Delta$, $F\on \nabla$ is 1-1, then 
$ F[\Delta] \cap F[\nabla ] = \emptyset $, or $F$ is symmetrical
($F(x,y)=F(y,x)$). 
\end{enumerate}
\end{Lemma}

\begin{proof}  1 and 2 are easy.   For 3, assume that $F(\alpha ,
\beta ) = F(\delta , \gamma )$, with $\alpha < \beta$, $\gamma <
\delta $. We have to distinguish several cases: 
\begin{description}
\item [Case 1] $\alpha = \gamma$, $\beta = \delta $.  Since
$F(\alpha,\beta)=F(\beta,\alpha)$, and $F$ is canonical, we have
$F(x,y)=F(y,x)$ for all $x,y$, so $F$ is symmetrical. 
\item [Case 2] $ \alpha = \gamma < \beta < \delta $. So
$F(\alpha,\beta)= F(\delta,\alpha)$. Pick any  $\beta ' , \delta '
$ with  $\delta < \beta ' < \delta ' $. 
\\ Then $(\alpha,\beta,\gamma,\delta) \sim (\alpha, \beta, \gamma,
\delta ' )$, 
so  $F(\alpha,\beta)=F(\delta ,\alpha)$ implies 
$ F(\alpha,\beta) = F(\delta ' , \alpha)$, this means 
$F(\delta ' , \alpha) = F (\delta , \alpha)$. 
\\
Similarly we find 
$ F(\alpha,\beta)  = F(\alpha,\beta')$.  So $F$ is neither 1-1 on
$\Delta$ nor 1-1 on $ \nabla $. 
\item [Other cases]  Similar to case 2.  
\end{description}
\end{proof}

\begin{Lemma}\label{lemma.HH}
Let $\C$ be a clone containing all unary functions.  If $\C$ contains
a heavily binary function $H$ and also a canonical function $F$ which
is 1-1 on ${\Delta}$, 
 then~$\C=\O$. 
\end{Lemma}

\begin{proof}
By \ref{fact.zorn}, 
it is enough to find a function $g\in \C$ which is 1-1 on~$\nd$.

If $F$ is symmetrical and 1-1 on $\Delta$ (and also 1-1 on~$\nabla$,
of course), then we may assume
(replacing $F$ by $h\circ F$ for some appropriate $h\in \oo 1$, if
necessary), that 
$F(x,y) > \max(x,y)$ for all~$x,y$. We claim that  the function
$$ (x,y) \mapsto F(x, F(x,y))$$
is 1-1 on~$\nd$.   Indeed, if $F(x,F(x,y))=F(x', F(x',y'))$, then we
have: 
\\ either~$x=x'$, $F(x,y)=F(x',y')$,
\\ or $x=F(x',y')$, $F(x,y)=x'$. 
\\ 
In the first case we get either  $y=y'$ directly, or~$x=y'$, $y=x$, so
again~$y=y'$.
\\
The second case leads to a contradiction:  $x=F(x',y')> x'$, $x <
F(x,y)=x'$.  

So we assume now that $F: \lambda \times \lambda \to \lambda $
 is canonical but not symmetrical.   By lemma~\ref{lemma.canonical},
 we know that $F[\Delta ] \cap F[\nabla ] = \emptyset$.  Replacing $F$
 by $h\circ F$ for an appropriate $h\in \oo 1$, we may assume that 
\begin{itemize}
\item $F\on \Delta$ is constantly~$0$.
\item $F\on \nabla $ takes only even values~$>0$, and is 1-1. 
\end{itemize}
Since $\C$ contains a heavily binary function, $\C$ contains some
function $H$ with $H(0,x)=x=H(x,1)$ for all~$x>0$.
 Now check that the map 
$(x,y)\mapsto H( Fxy, Fyx + 1 ) $ is a pairing function. 
\end{proof}

\subsection*{Proof of the theorem}

Assume that $\tau$ is a term for a function in $cl(\C_1\cup \C_2)$
  representing a 1-1 function on $\lambda
\times \lambda $.  Find a set $ S \subseteq \lambda $ of size $\lambda
$ such that  $\tau\on S$ is canonical (see definition
\ref{def.canonical}).  Since $\C$ contains all unary functions, $\C$
also contains a monotone 
 bijection between $S$ and $\lambda$, so wlog we will
assume that~$\tau$, as well as every subterm of~$\tau$, 
 is canonical on~$\lambda $. 

Let $ \Theta $ be the set of subterms of~$\tau$.

Let $ U_\Delta \subseteq \Theta  $ (and $U_\nabla \subseteq \Theta $)
 be the set of those terms $ \sigma $ which induce unary functions 
on $ \Delta$ ($ \nabla$, respectively), i.e., 
$$ U_\Delta = \{ \sigma\in \Theta:  \exists f\in  \lambda^\lambda  , 
\ 
 \left[ \forall (\alpha,\beta)\in \Delta:
 \sigma(\alpha,\beta)=f(\alpha)\right] 
\mbox { or } 
 \left[ \forall (\alpha,\beta)\in \Delta:
 \sigma(\alpha,\beta)=f(\beta )\right] \} $$
Let $ \sigma $ be a minimal subterm of $ \Theta $
 which is not in $ U_\Delta \cap
U_\nabla$, wlog~$ \sigma \notin U_\nabla$. 

Let $ G$ be the outermost function in the term $ \sigma$, say $ G\in
\C_1$, $ G $ $n$-ary.
   It remains to show that $\C_1 $ contains a pairing
function. 

All  proper subterms of $ \sigma$ represent unary functions, so there
are 
$n$ functions $ f_1, \ldots, f_n$ and some $k \le n$ with 
$$ \forall \alpha < \beta: \ \sigma(\alpha,\beta) = 
G(f_1(\alpha), \ldots, f_{k-1}(\alpha), f_k(\beta),\ldots,
f_{n}(\beta)).$$

So the function induced by $\sigma$ (which we again call $\sigma$) is
in~$\C_1$.  Now $\sigma\on \Delta$ is not essentially  unary. But $\sigma$
is canonical, so by lemma \ref{lemma.HH} we have a pairing
function in~$\C_1$. 

\section{Appendix: set theoretic assumptions} \label{section.app}

\begin{Definition}  Let $\lambda $, $\mu$, $n$, $c$ be cardinals (usually:  $\lambda $ and $\mu$ infinite, $n$ finite). 
The ``partition symbol''
$$  \lambda \to (\mu)^ n_c$$
says:  Whenever the set $[\lambda]^ n$, the set of subsets of $\lambda
$ of cardinality $n$ is partitioned into $c$ classes (i.e., whenever
$f:[\lambda ]^ n\to C$, where $|C|=c$), then there is a subset $A
\subseteq \lambda $ with at least $\mu $ many elements such that  all
subsets of $A$ of size $n$ are in the same equivalence class (i.e.,
the restriction of $f$ to $[A]^ n$ is a constant function).  
\end{Definition}

For example, the infinitary Ramsey theorem 
$$  {\aleph_0} \to ({\aleph_0} )^2_2$$
says: whenever the edges of a complete (undirected)
 graph on countably many
vertices are colored with $2$ colors, then there is an infinite
complete subgraph, all of whose edges have  the same color.

We will mainly be interested in the situation~$\lambda \to (\lambda)^ 2_2$.
If $\lambda \to (\lambda)^2_2$, and $\lambda$ is an uncountable
cardinal, then $\lambda $ is called ``weakly compact''. 

\begin{Fact} \label{fact.allfinite}
If $\lambda \to (\lambda)^ 2_2$, then for all finite $n,c$ we have 
$\lambda \to (\lambda)^ n_c$.  
\end{Fact}
(In  fact \ref{fact.reader}, we use this property
in the particular case $n=4$ and some large number
$c$, approximately  $c=3^{256}$.)

The property $\lambda \to (\lambda )^ 2_2$ is a rather strong
statement, i.e., it has many interesting consequences.    Therefore,
its  mere negation, 
$$ \lambda \not \to (\lambda )^ 2_2$$
or explicitly: 
\begin{quote} There is a map $f:[\lambda ]^2 \to \{0,1\}$ such that 
for any $A \subseteq \lambda $ of cardinality $\lambda$ the function
$f\on [A]^ 2 $ is not constant [i.e., is onto $\{0,1\}$]
\end{quote}

is a rather weak property of~$\lambda$.   There is, however, a
strengthening of this negative partition relation which already
yields interesting consequences. 

\begin{Definition}
The statement $\lambda \not\to[\lambda]^ 2_\lambda  $, the ``negative square
bracket partition relation''  means:
\begin{quote} There is a map $f:[\lambda ]^2 \to \lambda $ such that 
for any $A \subseteq \lambda $ of cardinality $\lambda$ the function
$f\on [A]^ 2 $ is {\em onto} $\lambda$
\end{quote}
\end{Definition}

We will now consider an even stronger property of $\lambda$:

\begin{Definition}\label{def.pr}
Let $\lambda\ge {\aleph_0}  $
 and $\mu$ be cardinals. The statement (or,
depending on your point of view, the ``principle'' or ``axiom'' 
 $Pr(\lambda,\mu)$ is defined as follows: 
\begin{quote}
There is a symmetric  $c:\lambda \times \lambda \to \mu$ with the
following property: 
\\
For all~$k\in \omega$, for all sequences $(a_i:i< \lambda )$ of
pairwise disjoint   subsets of $\lambda $,
of size~$k$, 
for all $c_0\in \mu$:
$$ \mbox{ there are $i<j<\lambda $ such that  
$c \on (a_i\times a_j) $ is constant with value $c_0$}$$
\end{quote}
\end{Definition}

Note that if we consider the case $\mu=\lambda$, and weaken the
conclusion by allowing only $k=1$, we get just $\lambda
\not\to[\lambda]^2_\lambda $.

This statement as well as several variants of it, are discussed in 
\cite[III.4 and appendix 1]{Sh:g}.   What we call $Pr(\lambda,\mu)$
corresponds to $Pr_1(\lambda,\lambda,\mu, {\aleph_0})$ there.

\bigskip

While the property $Pr(\lambda,\lambda)$ is quite strong 
(in particular:  sufficiently strong  to prove the result in 
section~\ref{section.many}), it turns out that is is not so rare: 
$Pr(\lambda, \lambda)$  holds for many successor
cardinals already in ZFC without extra axioms. 
More general results (with proofs) 
 can be found in chapter III of \cite{Sh:g}, and also 
in \cite{Sh:535} and \cite{Sh:572}. 

\begin{enumerate} 
\item If there is a nonreflecting $S \subseteq \{\delta < {\aleph_2} :
cf(\delta) = {\aleph_0} \} $, then $Pr(\lambda,\lambda)$ holds. 
See  \cite[III.4.6C(6)]{Sh:g}.


\item If $ \kappa \ge {\aleph_1} $ is regular, then~$Pr(\kappa^+, \kappa^+)$. 
See  \cite[III.4.8(1)]{Sh:g}, and \cite[theorem 1.1]{Sh:572} for the proof 
of~$Pr(\aleph_2,\aleph_2)$. 

\item  If $\kappa $ is singular, and the set of Jonsson cardinals
(=cardinals without a Jonsson algebra) is
bounded in  $\kappa$, then $Pr(\kappa^ +, \kappa ^ +)$ holds.  In
particular, $Pr(\aleph_{\omega+1},\aleph_{\omega+1})$ holds. 
See \cite[1.18]{Sh:535}. 

\end{enumerate}

\ifdraft
\bibliographystyle{lit-unsrt}
\else 
\bibliographystyle{plain}
\fi

\bibliography{listb,listx,lista}

\end{document}